\documentclass{aptpub}

\authornames{ANDREAS E. KYPRIANOU AND XIAOWEN ZHOU} 
\shorttitle{General tax structures and the L\'evy insurance risk model} 

\begin{document}

\title{GENERAL TAX STRUCTURES AND \\THE L\'EVY INSURANCE RISK MODEL.} 

\bigskip

\authorone[The University of Bath]{Andreas E. Kyprianou} 
\addressone{Department of Mathematical Sciences,
The University of Bath,
Claverton Down,
Bath BA2 7AY, UK.
email: a.kyprianou@bath.ac.uk}
\authortwo[Concordia University]{Xiaowen Zhou} 
\addresstwo{Department of Mathematics and Statistics, Concordia University, 1455 de Maisonneuve Blvd W.,  Montr\'eal Qu\'ebec, H3G 1M8, Canada. email: xzhou@mathstat.concordia.ca}

\begin{abstract}
In the spirit of \cite{AH, ARZ} we consider a L\'evy insurance risk
model with tax payments of a more general structure than in the
aforementioned papers that was also considered in \cite{ABBR}. In terms of scale functions, we establish three fundamental
identities of interest which have stimulated a large volume of
actuarial research in recent years. That is to say, the two sided
exit problem,  the net present value of tax paid until ruin as well
as a generalized version of the Gerber-Shiu function. The method we
appeal to differs from \cite{AH, ARZ} in that we appeal
predominantly to excursion theory.
\end{abstract}

\keywords{Reflected L\'evy processes, passage problems, integrated exponential L\'evy processes, insurance risk processes, ruin, excursion theory.} 

\ams{60K05, 60K15, 91B30}{60G70, 60J55} 

\section{Introduction and main results}
Recent advances in the analysis of the ubiquitous  ruin problem from the theory of insurance risk has seen a tendency to replace the classical Cram\'er-Lundberg surplus process with a general spectrally negative L\'evy processes; see for example \cite{ARZ, HPSV, KKM, RZ} to name but a few. In that case the surplus process is commonly referred to as a {\it L\'evy insurance risk process}. Although moving to this more complex setting, arguably,  does not  bring any more realistic features to the table than are already on offer in the classical Cram\'er-Lundberg model, a clear mathematical advantage has emerged. Working with L\'evy insurance risk processes forces one to approach the problem of ruin via excursion or fluctuation theory which does not use specific features of the underlying L\'evy process other than a generic path decomposition of the process in terms of excursions from its maximum which manifests itself in the form of a Poisson point process.

In this paper, we continue in this vein and build on ideas of L\'evy insurance risk processes with tax which were introduced and studied in \cite{AH, ARZ, ABBR}. Specifically we introduce a more general tax structure and therewith we establish, for the aggregate surplus process, new identities for the two sided exit problem, a generalized version of the Gerber-Shiu function as well the net present value of tax paid until ruin.

Henceforth the process $X=\{X_t: t\geq 0\}$ with probabilities $\{\mathbb{P}_x : x\in\mathbb{R}\}$ and natural filtration $\{\mathcal{F}_t : t\geq 0\}$ will denote a spectrally negative L\'evy process with the usual exclusion of processes in the latter class which have monotone paths (that is to say a pure increasing linear drift and the negative of a subordinator).  For convenience we shall always denote $\mathbb{P}_0$ by $\mathbb{P}$.
Let
\[
\psi(\lambda) = \log\mathbb{E}(e^{\lambda X_1})
\]
be the Laplace exponent of $X$ which is known to be finite for at least $\theta \in[0,\infty)$ in which case it is a strictly convex and infinitely differentiable function. The asymptotic behaviour of $X$ is characterized by $\psi'(0+)$,
so that $X$ drifts to $\pm\infty$ (oscillates) accordingly as $\pm\psi'(0+)>0$ ($\psi'(0+)=0$). When $X$ plays the role of the surplus process, it is usual to make the assumption that $\psi'(0+)>0$ which is equivalent to the {\it net profit condition} in the case that $X$ is a Cram\'er-Lundberg process. However, this condition is not necessary for any of the forthcoming analysis.

Denote by $S=\{S_t : t\geq 0\}$ the process which describes the running supremum of $X$, that is to say, $S_t = \sup_{s\leq t}X_s$ for
 each $t\geq 0$. Following \cite{ABBR}, we are interested in modeling tax payments from the L\'evy insurance process in such a way that the cumulative payment until
  time $t$ is given by
\[
\int_0^t \gamma(S_u)dS_u
\]
where $\gamma:[0,\infty)\to [0,1)$  is a measurable function which satisfies
\begin{equation}
\int^\infty_0 (1 - \gamma(s))ds =\infty
\label{fortheinverse}
\end{equation}
In that case the aggregate surplus process, the primary object of our study, is given by
 \begin{equation}
 U_t:=X_t-\int_0^t \gamma(S_u)dS_u.
 \label{agg}
 \end{equation}
In the special case that $\gamma$ is a constant in $(0,1)$ our L\'evy insurance risk process with tax agrees with the model introduced in \cite{AH, ARZ}.  In the case that $\gamma=0$, we are back to a regular L\'evy insurance risk process.  Note also that processes of the form (\ref{agg}) constitute a subclass of controlled L\'evy risk processes, the latter being of popular interest in recent literature; see for example \cite{APP, KL, KP}

 \bigskip

 In order to state the results alluded to above which concern path functionals of $U$, we must first introduce more notation.  As is now usual when studying L\'evy risk processes, a key element of the analysis involves the use of scale functions, defined as follows.
For every $q\geq 0$ there exists a  function $W^{(q)}:\mathbb{R}\rightarrow [0,\infty)$ such that $W^{(q)}(x)=0$ for all $x<0$ and otherwise is almost everywhere differentiable
on $[0,\infty)$ satisfying,
\begin{equation}\label{eq:scale}
\int_0^\infty e^{-\lambda x}
W^{(q)}(x)dx\,=\,\frac{1}{\psi(\lambda)-q},\text{ { } for { }}
\lambda>\Phi(q),
\end{equation}
where $\Phi(q)$ is the largest solution to the equation
$\psi(\theta)=q$ (there are at most two). We shall write for short
$W^{(0)}=W$ . It is known that when $X$ has paths of unbounded variation, the scale functions $W^{(q)}$ are continuously differentiable on $(0,\infty)$, and when $X$ has paths of bounded variation, they are almost everywhere differentiable. In either case we shall denote by $W^{(q)\prime}$ the associated density. It is also known that if $X$ has a Gaussian component then $W^{(q)}$ is twice continuously differentiable on $(0,\infty)$.

There exists a well known exponential change of measure that one may perform for spectrally negative L\'evy processes,
\begin{equation}
\left.\frac{d\mathbb{P}_x^\vartheta}{d\mathbb{P}_x}\right|_{\mathcal{F}_t}
= e^{\vartheta (X_t-x)- \psi(\vartheta)t}
\label{COM}
\end{equation}
for $x\in\mathbb{R}$ and $\vartheta\geq 0$ under which $X$ remains within the class of
spectrally negative L\'evy processes. In particular if $\nu(dx)$ is the L\'evy measure of $-X$ under $\mathbb{P}$ then $e^{-\vartheta x}\nu(dx)$ is its L\'evy measure under $\mathbb{P}^\vartheta$.
It will turn out  to be useful to introduce an additional parameter to the scale
functions described above in the light of this change of measure.
Henceforth we shall refer to the functions $W_\vartheta$  where $\vartheta\geq 0$ as the
functions that play the role of the scale functions defined in the
previous paragraph but when considered under the measures
$\mathbb{P}^\vartheta$.

Next define
\[
\tau_a^+ := \inf \{t > 0 \colon U_t > a \} \text{ and }\tau_0^-
:= \inf \{t
> 0 \colon U_t < 0 \}
\] with the convention
$\inf\emptyset=\infty$.  For $s\geq x$ define
 \begin{equation}
 \bar{\gamma}(s):=s-\int_x^s\gamma(y)dy=x+\int_x^s(1-\gamma(y))dy.
 \label{gamma-bar}
 \end{equation}
By differentiating (\ref{gamma-bar}) latter we note that since $\gamma\in[0,1)$ it follows that $\bar\gamma$ is a strictly increasing function. Moreover, since it is continuous it has a well defined inverse on $[x,\infty)$ which we denote by $\bar\gamma^{-1}$.

We may now present the three main results of this paper as promised. Their proofs will be given in the subsequent sections.

For each $x>0$, the process $L_t : = S_t-x$, $t\geq 0$, serves as a local time at $0$ for the Markov process $ Y: = S-X$ under $\mathbb{P}_x$. Write
$L^{-1}:=\{L^{-1}_t:t\geq 0\}$ for the right continuous inverse of $L$.

\begin{thm}[Two sided exit problem]\label{T:survival}
For any $0 < x < a$,  we have
\begin{equation}\label{E:laplace}
\bE_x\left[ e^{-q \tau^+_a} \ind_{\{\tau_a^+ < \tau_0^-\}} \right]
= \exp\left\{-\int_x^a
\frac{W^{(q)\prime}(y)}{W^{(q)}(y)(1-\gamma(\bar{\gamma}^{-1}(y)))}dy \right\}.
\end{equation}
\end{thm}

\bigskip

\begin{thm}[Net present value of tax paid until ruin]\label{NPV} For any $0 < x < a$,  we have
\begin{equation}\label{dividends}
\bE_x\left[
\int_0^{\tau^-_0} e^{-qu} \gamma(S_u) dS_u \right]= \int_x^\infty \exp\left\{- \int_x^t \frac{W^{(q)\prime}(\bar\gamma(s))}{W^{(q)} (\bar\gamma(s))} ds\right\} \gamma(t)dt.
\end{equation}
\end{thm}

\bigskip

\begin{thm}\label{G-S} For each $t\geq 0$ let  $S_t^U:=\sup_{s\leq t}U_s$ and
let $\kappa = L^{-1}_{L_{\tau^-_0}-}$, the last moment that tax is paid before ruin.
Denote by $\nu$  the jump measure of $-X$.
For any $y,z>0, \theta>y$ and $\alpha,\beta\geq 0$, we have
\begin{eqnarray}\label{gs}
\lefteqn{
\mathbb{E}_x\left(
e^{-\alpha\kappa - \beta(\tau^-_0 - \kappa)}; S^U_{\tau^-_0}\in d\theta, U_{\tau^-_0-}\in dy, - U_{\tau^-_0}\in dz
\right)
}&&\notag\\
&&=
\frac{1}{1- \gamma(\bar\gamma^{-1}(\theta))}  \exp\left\{-\int_x^\theta
\frac{W^{(\alpha)\prime}(y)}{W^{(q)}(\alpha)(1-\gamma(\bar{\gamma}^{-1}(y)))}dy \right\} \notag\\
&&\hspace{3cm}\left\{W^{(\beta)\prime}(\theta-y)-\frac{W^{(\beta)\prime}(\theta)}{W^{(\beta)}(\theta)}W^{(\beta)}(\theta-y)\right\}\nu(y+dz)d\theta dy.
\end{eqnarray}
Furthermore, we also have,
\begin{eqnarray}
\label{gs-creep}
\lefteqn{\mathbb{E}_x\left(
e^{-\alpha\kappa - \beta(\tau^-_0 - \kappa)}; S^U_{\tau^-_0}\in d\theta, U_{\tau^-_0} =0\right) }\notag\\
&&=
\frac{1}{1- \gamma(\bar\gamma^{-1}(\theta))}  \exp\left\{-\int_x^\theta
\frac{W^{(\alpha)\prime}(y)}{W^{(\alpha)}(y)(1-\gamma(\bar{\gamma}^{-1}(y)))}dy \right\}\notag\\
&&\hspace{5cm}\cdot\frac{\sigma^2}{2} \left\{ \frac{W^{(\beta)\prime}(\theta)^2}{W^{(\beta)}(\theta)} - W^{(\beta)\prime\prime}(\theta)\right\},
\end{eqnarray}
where $\sigma$ is the Gaussian coefficient in the L\'evy-It\^o decomposition of $X$.
\end{thm}

\begin{remark}\rm
When $\gamma\in(0,1)$ is a constant, we note that the expressions (\ref{E:laplace}) and (\ref{dividends}) agree with formulas (3.1) and (3.2) in \cite{ARZ}. Indeed we have $\bar\gamma(s)  = s(1-\gamma) + \gamma x$. For Theorem \ref{T:survival} we have
\[\bE_x\left[ e^{-q \tau^+_a} \ind_{\{\tau_a^+ < \tau_0^-\}} \right]
= \exp\left\{-\int_x^a
\frac{W^{(q)\prime}(y)}{W^{(q)}(y)(1-\gamma)}dy \right\}=\left(\frac{W^{(q)}(x)}{W^{(q)}(a)}\right)^{1/(1-\gamma)}.\]

For Theorem \ref{NPV} we have by two changes of variables
\begin{eqnarray*}
\int_x^\infty \exp\left\{- \int_x^t \frac{W^{(q)\prime}(\bar\gamma(s))}{W^{(q)} (\bar\gamma(s))} ds\right\} \gamma(t)dt &=&\gamma\int_x^\infty \exp\left\{- \frac{1}{1-\gamma}\int_x^{\bar\gamma(t)} \frac{W^{(q)\prime}(y)}{W^{(q)} (y)} dy\right\} dt \\
&=&\gamma\int_x^\infty \left(\frac{W^{(q)}(x)}{W^{(q)}(\bar\gamma(t))}\right)^{1/(1-\gamma)}dt\\
&=&\frac{\gamma}{1-\gamma}\int_x^\infty \left(\frac{W^{(q)}(x)}{W^{(q)}(u)}\right)^{1/(1-\gamma)}du,
\end{eqnarray*}
which is formula (3.2) in \cite{ARZ}.

Theorem \ref{G-S} on the other hand gives a new result for the setting of \cite{ARZ}. In particular, we have
\begin{eqnarray*}
\lefteqn{
\mathbb{E}_x\left(
e^{-\alpha\kappa - \beta(\tau^-_0 - \kappa)}; S^U_{\tau^-_0}\in d\theta, U_{\tau^-_0-}\in dy, - U_{\tau^-_0}\in dz
\right)
}&&\notag\\
&&= \frac{1}{1- \gamma}  \left(\frac{W^{(\alpha)}
(x)}{W^{(\alpha)}(\theta)}\right)^{1/(1-\gamma)}\left\{W^{(\beta)\prime}(\theta-y)-\frac{W^{(\beta)\prime}(\theta)}{W^{(\beta)}(\theta)}W^{(\beta)}(\theta-y)\right\}\nu(y+dz)d\theta
dy
\end{eqnarray*}
and
\begin{eqnarray*}
\lefteqn{\mathbb{E}\left(
e^{-\alpha\kappa - \beta(\tau^-_0 - \kappa)}; S^U_{\tau^-_0}\in d\theta, U_{\tau^-_0} =0\right) }\notag\\
&&= \frac{\sigma^2}{2(1-\gamma)}\left(\frac{W^{(\alpha)}
(x)}{W^{(\alpha)}(\theta)}\right)^{1/(1-\gamma)} \left\{
\frac{W^{(\beta)\prime}(\theta)^2}{W^{(\beta)}(\theta)} -
W^{(\beta)\prime\prime}(\theta)\right\}.
\end{eqnarray*}

\bigskip

Finally note that when $\gamma =0$ and the process $U$ agrees with the  L\'evy insurance risk process $X$, the last two formulae above give us two new expressions for the time value of the overall maximal wealth accumulated prior to ruin, the wealth immediately before ruin and the deficit at ruin.
\end{remark}

\begin{remark}\rm
One major criticism of working with scale functions is that, in principle, one has only solved the problems of interest up to inverting the Laplace transform (\ref{eq:scale}). However, in the last year there have been a number of developments in the theory of scale functions which has seen a large number of explicit examples appearing in the literature; including the case of Cram\'er-Lundberg models. See for example \cite{CKP, HK, KR, P}. The paper \cite{S} also gives recipes for evaluating scale functions numerically.
\end{remark}

\section{Proofs of Main results}
We begin this section by pointing out some important features of the running supremum of the aggregate process (\ref{agg}) which turns out to be key in our use of excursion theory in the forthcoming proofs.
\begin{lemma}\label{suprema}
We have that
\begin{equation}
S^U_t = S_t - \int_0^t \gamma(S_s)dS_s
\label{U-sup}
\end{equation}
and that
the random times $\{t\geq 0 : U_t = S^U_t\}$ agree precisely with $\{t\geq 0 : X_t = S_t\}$.
\end{lemma}
\begin{proof}
Note that on the one hand
\begin{equation}
S^U_t = \sup_{s\leq t} \left\{X_s - \int_0^s \gamma(S_u)dS_u\right\}\geq  \sup_{s\leq t} X_s - \int_0^t \gamma(S_u)dS_u = S_t - \int_0^t \gamma(S_u)dS_u.
\label{lower}
\end{equation}
On the other hand since $X_t\leq S_t$ we have
\[
U_t\leq  S_t -  \int_0^t \gamma(S_u)dS_u  = \int_0^t (1- \gamma(S_u))dS_u
\]
and hence, since $\gamma (y)\in[0,1)$ for all $y\geq 0$,
\begin{equation}
S^U_t \leq  \sup_{s\leq t} \int_0^s (1- \gamma(S_u))dS_u = \int_0^t (1- \gamma(S_u))dS_u = S_t  - \int_0^t \gamma(S_u)dS_u
\label{upper}
\end{equation}
Together (\ref{upper}) and (\ref{lower}) imply (\ref{U-sup}). Now suppose that $t'\in\{t\geq 0 : X_t = S_t\}$. This implies that
\[
U_{t'} = X_{t'} - \int_0^{t'}\gamma(S_u)dS_u=S_{t'} - \int_0^{t'}\gamma(S_u)dS_u = S^U_{t'}
\]
and hence $t'\in\{t\geq 0 : U_t = S^U_t\}$. On the other hand, if
$t''\in \{t\geq 0 : U_t = S^U_t\}$, then
\[
X_{t''} - \int_0^{t''}\gamma(S_u)dS_u  = U_{t''} = S^U_{t''} = S_{t''} - \int_0^{t''}\gamma(S_u)dS_u
\]
showing that $X_{t''} = S_{t''}$ and hence $t''\in\{t\geq 0: X_t = S_t\}$.
\hfill$\square$\end{proof}

For the remaining proofs we shall also make heavy use of excursion theory for the process $S-X$ for which
we refer to \cite{be96} for background reading.  We shall spend a moment here setting up some necessary notation
which will be used throughout the remainder of the paper.  The Poisson process of excursions indexed by local time
shall be denoted by $\{(t, \epsilon_t): t\geq 0\}$ where
\[
\epsilon_t = \{\epsilon_t(s) := X_{L^{-1}_{t}} - X_{L^{-1}_{t-}+s}: 0< s\leq L^{-1}_{t} - L^{-1}_{t-} \}
\]
 whenever $\sigma(\epsilon_t):=L^{-1}_{t} - L^{-1}_{t-}>0$.  Accordingly we refer to a generic excursion as $\epsilon(\cdot)$
 (or just $\epsilon$ for short as appropriate) belonging to the space  $\mathcal{E}$ of canonical excursions.
 The intensity measure of the process $\{(t, \epsilon_t): t\geq 0\}$ is given by $dt\times dn$ where $n$ is a measure on the space of
 excursions (the excursion measure). An $n$-measurable
functional of the canonical excursion which will be of prime interest is $\bar\epsilon= \sup_{s\geq 0}\epsilon(s)$.
A useful formula for this functional that we shall make use of is the following (cf. \cite{kyp06})
\begin{equation}
n(\overline{\epsilon}> x) = \frac{W'(x)}{W(x)}
\label{excursion-tail}
\end{equation}
providing that $x$ is not a point of discontinuity in the derivative of $W$ (which is only a concern when $X$ has paths of bounded variation, in which case there are at most a countable number).

Lemma \ref{suprema} also has an important bearing on the process of
excursions described above. Indeed,  from the identity (\ref{U-sup})
we note that if $L_t = s$, or equivalently $S_t = x+s$, under
$\mathbb{P}_x$, then  $S^U_t = \bar\gamma(x+s)$. Moreover,
$L^{-1}_{s} = \tau^+_{\bar\gamma(x+s)}$, or equivalently $\tau^+_a =
L^{-1}_{\bar\gamma^{-1}(a) -x}$, under $\mathbb{P}_x$ and the
excursions of $U$ away from its maximum agree precisely with $\{(t,
\epsilon_t): t\geq 0\}$.

\begin{proof}[Proof of Theorem \ref{T:survival}]
Taking account of the remarks following Lemma \ref{suprema} we have that
the
event $\{\tau_a^+ < \tau_0^-\}$ is the same as
\[
\{\bar\epsilon_s \leq  \bar{\gamma}(x+s) , \forall \, 0 \leq s <
\bar{\gamma}^{-1}(a)-x \}.
\]
 Then, for $x>0$,
\begin{equation*}
\begin{split}
\bP_x (\tau_a^+ < \tau_0^-) 
&=\mathbb{P}_x( \bar\epsilon_s \leq  \bar{\gamma}(x+s) , \forall \, 0 \leq
s<\bar{\gamma}^{-1}(a)-x )\\
&=\exp\left\{-\int_0^{\bar{\gamma}^{-1}(a)-x}   n(\bar\epsilon >\bar\gamma(x+s))ds\right\}\\
& = \exp \left\{-\int_0^{\bar{\gamma}^{-1}(a)-x} \frac{W^{\prime}(\bar{\gamma}(x+s))}{W(\bar{\gamma}(x+s))} \, ds \right\} \\
&=\exp \left\{-\int_x^{a}
\frac{W^{\prime}(y)}{W(y)(1-\gamma(\bar{\gamma}^{-1}(y)))} \, dy
\right\},
\end{split}
\end{equation*}
where we change the variable making use of the fact that, since
$\bar\gamma(\bar\gamma^{-1}(s))=s$, we have from the chain rule
\[
\frac{d}{ds} \bar\gamma^{-1}(s) =
\frac{1}{\bar\gamma'(\bar\gamma^{-1}(s))} =
\frac{1}{1-\gamma(\bar\gamma^{-1} (s))}.
\]
Next, note that
\begin{eqnarray}
\bP_x^{\Phi(q)}( \tau_a^+ < \tau_0^- ) &=&
\bE_x \left[ e^{\Phi(q)(X_{\tau^+_a} - x)-q \tau^+_a} \ind_{\{\tau_a^+ < \tau_0^-\}}\right]\notag\\
&=& \exp\left\{\Phi(q)\left(\bar{\gamma}^{-1}(a)-x\right)\right\}
\bE_x \left[ e^{-q \tau^+_a} \ind_{\{\tau_a^+ < \tau_0^-\}}\right]
\label{one}
\end{eqnarray}
where we have appealed to the change of  measure (\ref{COM}) with $\vartheta = \Phi(q)$ and the final equality follows by virtue of the fact that on $\{\tau^+_a <\infty\}$
\begin{eqnarray*}
 X_{\tau_a^+}&=&U_{\tau^+_a}+\int_0^{\tau_a^+} \gamma(S_u)dS_u\notag\\
&=& a+\int_0^{L^{-1}_{\bar\gamma^{-1}(a) -x}}\gamma(S_u)dS_u\notag\\
&=& a+\int_x^{\bar{\gamma}^{-1}(a)}\gamma(y)dy\notag\\
&=&\bar{\gamma}^{-1}(a)
\label{two}
\end{eqnarray*}
where in the second equality we have made the change of variable $y
= S^{-1}_u$.
Note also that it is known (cf. Chapter 8 of  \cite{kyp06}) that for $q,x\geq 0$,
\begin{equation}
W^{(q)}(x) = e^{\Phi(q)x}W_{\Phi(q)}(x)
\label{q-Phiq}
\end{equation}
and hence
\begin{equation}
\frac{W_{\Phi(q)}'(x)(x)}{W_{\Phi(q)}(x)}=\frac{{W^{(q)}}'(x)}{W^{(q)}(x)}-\Phi(q).
\label{three}
\end{equation}
Piecing together (\ref{one}), (\ref{two}) and (\ref{three}) we get
\begin{equation*}
\begin{split}
&\bE_x \left[ e^{-q \tau^+_a} \ind_{\{\tau_a^+ < \tau_0^-\}}\right]\\
&= \bP_x^{\Phi(q)} ( \tau_a^+ < \tau_0^- ) \exp\left\{-\Phi(q)\left(\bar{\gamma}^{-1}(a)-x \right)\right\} \\
&=\exp \left\{-\int_x^{a}
\frac{{W_{\Phi(q)}}^{\prime}(y)}{W_{\Phi(q)}(y)(1-\gamma(\bar{\gamma}^{-1}(y)))}
\, dy \right\}\exp\left\{-\Phi(q)\left(\bar{\gamma}^{-1}(a)-u \right)\right\}\\
&=\exp\left\{-\int_x^a\frac{W^{(q)\prime}(y)}{W^{(q)}(y)(1-\gamma(\bar{\gamma}^{-1}(y)))}dy\right\},
\end{split}
\end{equation*}
where we have also used the  fact that
\[\int_x^a \frac{1}{1-\gamma(\bar{\gamma}^{-1}(y))}dy=\bar{\gamma}^{-1}(a)-\bar{\gamma}^{-1}(x)=\bar{\gamma}^{-1}(a)-x.\]
The proof is now complete.
\hfill$\square$\end{proof}


\begin{proof}[Proof of Theorem \ref{NPV}]
The proof builds on the experience of the calculations in the previous proof. We note that  the process $S$ does not increase on the time interval $(L^{-1}_{L_{\tau^{-}_0}-}, \tau^-_0)$ and  hence
\begin{eqnarray*}
\lefteqn{\mathbb{E}_x\left[\int_0^{\tau^-_0} e^{-qu} \gamma(S_u)dS_u
\right]}\\
&&=\mathbb{E}_x\left[\int_0^{L^{-1}_{L_{\tau^-_0}-}} e^{-qu} \gamma(L_u+x)dL_u
\right]\\
&&=\mathbb{E}_x\left[\int_0^{\infty}\ind_{\{t< L_{\tau^-_0}\}} e^{-qL^{-1}_t} \gamma(t+x)dt
\right]\\
&&=\int_0^{\infty}   \mathbb{E}_x\left[  e^{-q L^{-1}_{t} } \ind_{\{
\bar\epsilon_s \leq  \bar{\gamma}(x+s) , \forall \, 0 \leq s \leq t \}}
\right]  \gamma(t+x)dt\\
&&=\int_0^{\infty}  e^{-\Phi(q)t}  \mathbb{P}^{\Phi(q)}_x( \bar\epsilon_s \leq  \bar{\gamma}(x+s) , \forall \, 0 \leq s \leq t )  \gamma(t+x)dt\\
&&=\int_0^{\infty}  e^{-\Phi(q)t}  \exp\left\{- \int_0^t n_{\Phi(q)}(\bar\epsilon > \bar\gamma(x+s)) ds\right\}   \gamma(t+x)dt\\
&&=\int_0^{\infty}  e^{-\Phi(q)t}  \exp\left\{- \int_0^t \frac{W_{\Phi(q)}'(\bar\gamma(x+s))}{W_{\Phi(q)}(\bar\gamma(x+s))} ds \right\}   \gamma(t+x)dt\\
&&=\int_0^{\infty}   \exp\left\{- \int_0^t \frac{W^{(q)\prime}(\bar\gamma(x+s))}{W^{(q)}(\bar\gamma(x+s))} ds\right\}   \gamma(t+x)dt
\end{eqnarray*}
where in the fifth equality the measure $n_{\Phi(q)}$  plays the role of $n$ under $\mathbb{P}^{\Phi(q)}$, in the penultimate equality we have used (\ref{excursion-tail}) and the final equality uses (\ref{three}). The proof is completed by applying a straightforward change of variables.
\hfill$\square$\end{proof}

Before turning the proof of Theorem \ref{G-S}, we need first to
prove an additional auxiliary result. To this end, define $\rho_a=
\inf\{s>0 : \epsilon(s)>a \}$, the first passage time above $a$ of
the canonical excursion $\epsilon$.
We also need the first passage times for the underlying L\'evy process $X$,
\[
T^+_x = \inf\{t>0 : X_t >x\}\text{ and }T^-_x=\inf\{t>0 : X_t <x\}
\]
for all $x\in\mathbb{R}$.

\begin{lem}\label{aux-lemma}
For any $y,z>0$ and $q\geq 0$, we have
\begin{eqnarray*}
\lefteqn{n\left(e^{-q\rho_a}; a- \epsilon(\rho_a-) \in dy,\epsilon(\rho_a ) -a\in dz \right) }&&\\
&&= \left\{W^{(q)\prime}(a-y)-\frac{W^{(q)\prime}(a)}{W^{(q)}(a)}W^{(q)}(a-y)\right\}\nu(y+dz)dy
\end{eqnarray*}
and
\[
n\left( e^{-q \rho_a}; \epsilon(\rho_a)=a \right)
=\frac{\sigma^2}{2} \left\{ \frac{W^{(q)\prime}(a)^2}{W^{(q)}(a)} -
W^{(q)\prime\prime}(a)\right\}.
\]
\end{lem}

\begin{proof} Recall that  $Y= S-X$. For the latter process introduce its first passage time
\[
\varsigma_a = \inf\{t>0 : Y_t >a\}.
\]
By a classical application of the compensation formula (see for example the treatment of a related problem in \cite{AKP}) we have for $q\geq 0$ that
\begin{eqnarray}
\lefteqn{\mathbb{E}(e^{-q \varsigma_a} ; a- Y_{\varsigma_a - } \in dy, Y_{\varsigma_a} - a\in dz)}&&\notag\\
&&=\mathbb{E}\left[\sum_{t\geq 0} e^{-q(L^{-1}_{t-}+ \rho_a(\epsilon_t))} \ind_{\{\sup_{s<t} \overline\epsilon_s \leq  a,   \rho_a(\epsilon_t) <\sigma(\epsilon_t) , a- \epsilon_t(\rho_a-) \in dy,\epsilon_t(\rho_a ) -a\in dz \} }
\right]\notag\\
&&=\int_0^\infty
\mathbb{E}\left[e^{-qL^{-1}_{t}} \ind_{\{ \sup_{s\leq t} \overline\epsilon_s \leq  a\}}\right] dt \cdot n\left(e^{-q\rho_a}; a- \epsilon(\rho_a-) \in dy,\epsilon(\rho_a ) -a\in dz \right)\notag\\
&&=\int_0^\infty e^{-\Phi(q)t} e^{-n_{\Phi(q)}(\bar\epsilon> a) t}dt \cdot n\left(e^{-q\rho_a}; a- \epsilon(\rho_a-) \in dy,\epsilon(\rho_a ) -a\in dz \right)\notag\\
&&=\int_0^\infty e^{-\Phi(q)t} \exp\left\{-\frac{W'_{\Phi(q) } (a)}{W_{\Phi(q)}(a) } t\right\}dt \cdot n\left(e^{-q\rho_a}; a- \epsilon(\rho_a-) \in dy,\epsilon(\rho_a ) -a\in dz \right)\notag\\
&&=\int_0^\infty \exp\left\{-\frac{W^{(q)\prime } (a)}{W^{(q)}(a) } t\right\}dt\cdot n\left(e^{-q\rho_a}; a- \epsilon(\rho_a-) \in dy,\epsilon(\rho_a ) -a\in dz \right)
\label{mine}
\end{eqnarray}
where in the first equality the time index runs over local times and the sum is the usual shorthand for integration with respect to the Poisson counting measure of excursions, and for the second equality we need the quasi left continuity for subordinator $L^{-1} $.

On the other hand, according to Theorem 1 of  \cite{Pist-Sem}, we have that
\begin{eqnarray}
\lefteqn{
\mathbb{E}(e^{-q \varsigma_a} ; a- Y_{\varsigma_a - } \in dy, Y_{\varsigma_a} - a\in dz)
}&&\notag\\
&&=\int_0^\infty \exp\left\{-\frac{W^{(q)\prime } (a)}{W^{(q)}(a) } t\right\}dt\notag\\
&&\hspace{1cm}\times \left(W^{(q)\prime}(a-y)-\frac{W^{(q)\prime}(a)}{W^{(q)}(a)}W^{(q)}(a-y)\right)\nu(y+dz)dy
\label{martijn's}
\end{eqnarray}
By comparing the left and right hand sides of (\ref{mine}) and (\ref{martijn's}) we thus have that
\begin{eqnarray*}
\lefteqn{n\left(e^{-q\rho_a}; a- \epsilon(\rho_a-) \in dy,\epsilon(\rho_a ) -a\in dz \right) }&&\\
&&= \left(W^{(q)\prime}(a-y)-\frac{W^{(q)\prime}(a)}{W^{(q)}(a)}W^{(q)}(a-y)\right)\nu(y+dz)dy
\end{eqnarray*}
as claimed.

For the proof of the second part we should first note that it known (cf. \cite{be96}) that the process $X$ creeps downwards
if and only if $\sigma\neq 0$. It thus follows that $Y$  creeps upwards if  and only if $\sigma\neq 0$.
Henceforth assume that $\sigma \neq 0$. The proof then follows the same reasoning except in (\ref{mine}) one replaces the
event $\{a- Y_{\varsigma_a - } \in dy, Y_{\varsigma_a} - a\in dz\}$ by $\{Y_{\varsigma_a} = a\}$ on the left hand side and $\{a- \epsilon(\rho_a-) \in dy,\epsilon(\rho_a ) -a\in dz\}$ by $\{\epsilon(\rho_a) =a\}$ on the right hand side. Furthermore as a replacement for  (\ref{martijn's}) in the argument we use instead
\[
\mathbb{E}(e^{-q\varsigma_a}; Y_{\varsigma_a} = a) = \frac{\sigma^2}{2} \left\{ \frac{W^{(q)\prime}(a)^2}{W^{(q)}(a)} - W^{(q)\prime\prime}(a)\right\}\cdot \int_0^\infty \exp\left\{-\frac{W^{(q)\prime } (a)}{W^{(q)}(a) } t\right\}dt,
\]
which is taken from Theorem 2 of \cite{Pist-Sem}.
\hfill$\square$\end{proof}

\begin{proof}[Proof of Theorem \ref{G-S}] We give the proof for the first identity. The proof of the second identity follows along exactly the same lines using the second part of Lemma \ref{aux-lemma} instead and is left as an exercise for the reader.

In a similar spirit to the proof of Lemma \ref{aux-lemma} we may write for a given open interval $B\subset(0,\infty)$,
\begin{eqnarray*}
\lefteqn{
\mathbb{E}\left(
e^{-\alpha\kappa - \beta(\tau^-_0 - \kappa)}; S^U_{\tau^-_0}\in B, U_{\tau^-_0-}\in dy, - U_{\tau^-_0}\in dz
\right)
}&&\\
&&=\mathbb{E}_x\left[
\sum_{t\geq 0}\ind_{\{L^{-1}_{t-} < \tau^-_0\}}e^{-\alpha L^{-1}_{t-} }
\ind_{\{S^U_{L^{-1}_{t-}}\in B\}}   e^{- \beta (\tau^-_0 - L^{-1}_{t-})}
\ind_{\{S^U_{\tau_0^-} = S^U_{L^{-1}_{t-}}\}}
\ind_{\{ U_{\tau^-_0 -} \in dy, - U_{\tau^-_0} \in dz\}}
\right].
\end{eqnarray*}
Note however that, on account of the fact that
\[
S^U_{L^{-1}_{t-}}=S^U_{L^{-1}_t} = x+t -
\int_x^{L^{-1}_t}\gamma(S_u)dS_u  = x+t -
\int_x^{x+t}\gamma(y)dy=\bar\gamma(x+t)
\]
we have that $\{S^U_{L^{-1}_t}\in B \} = \{\bar\gamma(x+t)\in B\}$.
Note also that $L^{-1}_{t} = \tau^+_{\bar\gamma(x+t)}$ and 
 $L^{-1}$ is quasi left continuous. Hence, applying the
compensation formula we have
\begin{eqnarray*}
\lefteqn{
\mathbb{E}\left(
e^{-\alpha\kappa - \beta(\tau^-_0 - \kappa)}; S^U_{\tau^-_0}\in B, U_{\tau^-_0-}\in dy, - U_{\tau^-_0}\in dz
\right)
}&&\\
&&=\mathbb{E}_x\Big[
\int_0^\infty  dt\cdot e^{-\alpha \tau^+_{\bar\gamma(x+t)} } \ind_{\{\bar\gamma(x+t)\in B\}} \ind_{\{\tau^+_{\bar\gamma(x+t)} <\tau^-_0 \}}\\
&&\hspace{1.5cm}\times n(e^{-\beta \rho_{\bar\gamma(x+t)}}; \bar\gamma(x+t) - \epsilon(\rho_{\bar\gamma(x+t)} -) \in dy, \epsilon(\rho_{\bar\gamma(x+t)})-\bar\gamma(x+t) \in dz)\Big]\\
&&=
\int_B  \frac{d\theta}{1- \gamma(\bar\gamma^{-1}(\theta))} \mathbb{E}_x\left[e^{-\alpha \tau^+_\theta }  \ind_{\{\tau^+_\theta <\tau^-_0 \}}\right]
 n(e^{-\beta \rho_\theta}; \theta - \epsilon(\rho_\theta -) \in dy, \epsilon(\rho_\theta)-\theta \in dz)
\end{eqnarray*}
where in the final equality we have applied a change of variable.

Now making use of the first part of Lemma \ref{aux-lemma} and the conclusion of Theorem \ref{T:survival} it thus follows that
\begin{eqnarray*}
\lefteqn{
\mathbb{E}\left(
e^{-\alpha\kappa - \beta(\tau^-_0 - \kappa)}; S^U_{\tau^-_0}\in d\theta, U_{\tau^-_0-}\in dy, - U_{\tau^-_0}\in dz
\right)
}&&\\
&&=
\frac{1}{1- \gamma(\bar\gamma^{-1}(\theta))}  \exp\left\{-\int_x^\theta
\frac{W^{(\alpha)\prime}(y)}{W^{(\alpha)}(y)(1-\gamma(\bar{\gamma}^{-1}(y)))}dy \right\} \\
&&\hspace{3cm}\left\{W^{(\beta)\prime}(\theta-y)-\frac{W^{(\beta)\prime}(\theta)}{W^{(\beta)}(\theta)}W^{(\beta)}(\theta-y)\right\}\nu(y+dz)d\theta dy
\end{eqnarray*}
as required.
\hfill$\square$\end{proof}

\subsection*{Acknowledgments}
The first author acknowledges the support of EPSRC grant number
EP/D045460/1. The second author is supported by an NSERC grant.

\end{document}